\documentclass[a4paper,12pt]{article}
\usepackage[psamsfonts]{amssymb}
\usepackage[psamsfonts]{eucal}

\usepackage{amsmath}
\usepackage{theorem}
\newtheorem{thm}{Theorem}
\newtheorem{prop}[thm]{Proposition}
\newtheorem{cor}[thm]{Corollary}
\newtheorem{lm}[thm]{Lemma}

\theorembodyfont{\rmfamily}

\theorembodyfont{\rmfamily}
\begin{document}

\title{Para-hyperhermitian surfaces}

\author{Johann Davidov\thanks{Partially supported by "L.Karavelov" Civil Engineering Higher School, Sofia, Bulgaria
under contract No 10/2009}, Gueo Grantcharov, Oleg
Mushkarov\thanks{Partially supported by CNRS-BAS joint research
project {\it Invariant metrics and complex geometry, 2008-2009}},
Miroslav Yotov}

\date{}
\maketitle

\centerline{\it Dedicated to Professor Stere Ianu\c{s}~on the
occasion of his 70th birthday}

 \rm
\begin{abstract}
In this note we discuss the problem of existence of
para-hyperhermitian structures on compact complex surfaces. We
construct examples of para-hypercomplex structures on Inoue surfaces
of type $S^{-}$ which do not admit compatible metrics.\\



\end{abstract}

\vspace{0.5cm}

\section{\bf Introduction}

Hypercomplex and hyperk\"ahler structures have been studied for a
long time and many interesting results and relations with other
fields have been established. Recently there is a growing interest
in their pseudo-Riemannian counterparts too due to the fact that
important geometry models of string theory carry such structures
\cite {OV}.
 The para-hyperhermitian structures arise as a
pseudo-Riemannian  analog of the hyperhermitian structures and it is
well known \cite{Kamada} that in four dimensions they lead to
self-dual metrics of neutral signature. There are many other
similarities between these two structures, but there are also
significant differences. For example, the para-hypercomplex
structures, the neutral analog of hypercomplex structures, exist in
any even dimension (not only in that divisible by 4) and, in
contrast to the latter, they may not have compatible metrics.

 An almost para-hypercomplex structure is a triple $(J_1,J_2,J_3)$ of anti-commuting
endomorphisms of the tangent bundle with $J_1^2=-J_2^2=-J_3^2=-Id$.
When the structures $J_1, J_2, J_3$ are integrable it is called
para-hypercomplex. There are two natural classes of metrics
compatible with such structures. The first one consists of the
neutral metrics for which the structure $J_1$ is an isometry while
$J_2$ and $J_3$ are anti-isometries. These metrics, called
para-hyperhermitian in this paper, give rise to three 2-forms
defined in the same way as the K\"ahler forms in the positive
definite case. When they are closed the structure is called
hypersymplectic \cite{Hit1990}, para-hyperk\"ahler \cite{Bl-V},
hyper-parak\"ahler \cite{kuche}, neutral hyperk\"ahler \cite{Ka},
pseudo-hyperk\"ahler \cite{DW}, etc. The second class consists of
positive definite metrics, for which the structures $J_1, J_2, J_3$
are isometries. Such metrics always exist and the analog of the
K\"ahler form for the structure $J_2$ is a symmetric form which is
in fact a neutral metric such that $J_1$ and $J_3$ are
anti-isometries but $J_2$ is an isometry. As pointed out in
\cite{cecilia}, neutral metrics with this property are interesting
in relation with the doubled geometry models of string theory,
introduced by C. Hull \cite{Hull}. We should note however that the
existence of para-hyperhermitian metrics leads to some additional
obstructions and a purpose of this note is to clarify the problem
for their existence.

 The paper is organized as follows. After the preliminary definitions (Section 2)
we recall in Section 3 Kamada's classification \cite{Ka} of compact
para- hyperk\"ahler surfaces and relate them to the existence of
parallel null vector fields. Then in Section 4 we establish some
necessary conditions for the existence of a para-hyperhermitian
metric
 with respect to a given para-hypercomplex structure on a 4-manifold and show that any
 two such metrics are conformally equivalent. In the last section we show that the Inoue
 surfaces of type $S^+$ have para-hyperhermitian structures and provide examples of
 para-hypercomlex structures  on Inoue surfaces of type $S^-$ which do not admit compatible
 para-hyperhermitian metrics.

\section{\bf Preliminaries}

 A pseudo-Riemannian metric on a smooth 4-manifold $M$ is called $neutral$ if it has signature
$(+,+,-,-)$. Unlike the Riemannian case, there are topological
restrictions for existence of a neutral metric on a compact manifold
since it is equivalent to existence of a field of tangent 2-planes
\cite{Steen}. We refer to \cite{Matsushita1991} for further
information in this direction.

\smallskip

  An almost para-hypercomplex structure on a smooth $4$-manifold $M$
consists of three endomorphisms $J_1,J_2,J_3$ of $TM$ satisfying the
relations
\begin{equation}\label{3.2}
J_1^2=-J_2^2=-J_3^2=-Id, \hspace{.1in} J_1J_2=-J_2J_1=J_3
\end{equation}
of the imaginary units of the paraquaternionic algebra (split quaternions). A
metric $g$ on $M$ is called compatible with the structure $\{J_1,J_2,J_3\}$ if
\begin{equation}\label{3.3}
g(J_1X,J_1Y)=-g(J_2X,J_2Y)=-g(J_3X,J_3Y)=g(X,Y)
\end{equation}
(such a metric is necessarily of split signature). In this case we
say that $\{g,J_1,J_2,J_3\}$ is an almost para-hyperhermitian
structure. For any such a structure we define three $2$-forms
$\Omega_i$ setting $\Omega_i(X,Y) = g(J_iX,Y)$, $i=1,2,3$. If the
Nijenhuis tensors of $J_1,J_2,J_3$ vanish, the structure
$\{g,J_1,J_2,J_3\}$ is called para-hyperhermitian. When additionally
the 2-forms $\Omega_i(X,Y) = g(J_iX,Y)$ are closed, the
para-hyperhermitian structure is called para-hyperk\"ahler. It is
well known \cite{Kamada} that the para-hyperhermitian metrics are
self-dual, whereas the para-hyperk\"ahler metrics are self-dual and
Ricci-flat.

\smallskip

  It is an observation of Hitchin \cite{Hit1990} (see also \cite{Ka})
that any para-hyperk\"ahler structure is uniquely determined by three
symplectic forms  $(\Omega_1, \Omega_2 , \Omega_3)$ satisfying the relations
$$
-\Omega_1^2=\Omega_2^2=\Omega_3^2, \quad \Omega_l\wedge\Omega_m = 0,\> l\neq m.
$$

A similar characterization holds for para-hyperhermitian structures
\cite{kuche,Ka}. They are uniquely determined by three
non-degenerate $2$-forms $(\Omega_1, \Omega_2 , \Omega_3)$ and a
$1$-form $\theta$ such that
\begin{equation}\label{phc}
-\Omega_1^2=\Omega_2^2=\Omega_3^2, \quad \Omega_l\wedge\Omega_m = 0,\> l\neq
m,\quad d\Omega_l=\theta\wedge\Omega_l.
\end{equation}

For any para-hyperhermitian structure on a $4$-manifold $M$, the $2$-form $\Omega=\Omega_2+i\Omega_3$ is of type
$(2,0)$ with respect to the complex structure $J_1$, hence the canonical bundle of the complex manifold
$(M,J_1)$ is smoothly trivial. Using the well-known classification of compact complex surfaces it follows that
para-hyperhermitian structures can exist only on the following surfaces: complex tori, K3 surfaces, primary
Kodaira surfaces, Hopf surfaces, Inoue surfaces of type $S_M, S_N^{\pm}$ and properly elliptic surfaces of odd
first Betti number. Note that except the $K3$ surfaces all these surfaces can be represented as quotients of Lie
groups factored by cocompact discrete subgroups (more details will appear in \cite{inprogress}).

\section{Para-hyperk\"ahler surfaces}

As is well known (c.f. \cite{bes}), any compact hyperk\"ahler
surface is either a complex torus with a flat metric or a
$K3$-surface with Calabi-Yau metric. In the neutral case, the
$(2,0)$-form $\Omega=\Omega_2+i\Omega_3$ is holomorphic (even
parallel), so the canonical bundle is holomorphically trivial. Using
this fact H.Kamada \cite{Ka} proved the following
\begin{thm}
If  $(M,g,J_1,J_2,J_3)$ is a compact para-hyperk\"ahler surface, then the
complex surface $(M,J_1)$ is biholomorphic to a complex torus or a primary
Kodaira surface.
\end{thm}

 Moreover, Kamada \cite{Ka, Kamada} obtained  a description of all para-hyperk\"ahler  structures on
both types of surfaces.

\begin{thm}
For any para-hyperk\"ahler structure on a complex torus $M={\Bbb
C}^2/\Gamma$ there are complex coordinates $(z_1, z_2)$ of ${\Bbb
C}^2$, such that the structure is defined by means of the following
symplectic forms:
$$\Omega_1 = Im(dz_1\wedge d\overline{z}_2)+(i/2)\partial\overline\partial\varphi, \\
$$
$$
\Omega_2 = Re(dz_1\wedge dz_2),\hspace{.1in} \Omega_3 = Im(dz_1\wedge dz_2),
$$
where $\varphi$ is a smooth function such that
\begin{equation}\label{kod-hk}
4i(Im(dz_1\wedge d\overline{z}_2)\wedge\partial\overline{\partial}\varphi =
\partial\overline{\partial}\varphi\wedge\partial\overline{\partial}\varphi.
\end{equation}
Conversely, any three forms $\Omega_1,\Omega_2,\Omega_3$ of the form
given above determine a para-hyperk\"ahler structure on the torus.
Moreover, its metric  is flat if and only if $\varphi$ is constant.
 \end{thm}

  Let us note that if $M$ is a product of two elliptic curves, then there are non-trivial
solutions of the equation (\ref{kod-hk}) (\cite{Kamada}) and it is
not known if such solutions exist when $M$ is not a product.

\smallskip

 Before stating Kamada's result about primary Kodaira surfaces, we recall their
definition.

 Consider the affine transformations
 $\rho_i(z_1,z_2) = (z_1+a_i,z_2+\overline{a}_iz_1+b_i)$ of ${\Bbb C}^2$,
 where $a_i$,$b _i$, $i=1,2,3,4$, are  complex
numbers such that $a_1=a_2=0, Im(a_3{\overline a}_4) = b_1$. Then
$\rho_i$ generate a group $G$ of affine transformations acting
freely and properly discontinuously on ${\Bbb C}^2$. The quotient
space $M={\Bbb C}^2/G$ is called a primary Kodaira surface.

\begin{thm}
For any para-hyperk\"ahler structure on a primary Kodaira surface
$M$ there are complex coordinates $(z_1, z_2)$ of ${\Bbb C}^2$ such
that the structure is defined by means of the following symplectic
forms:
$$\Omega_1 = Im(dz_1\wedge d\overline{z}_2)+iRe(z_1)dz_1\wedge d\overline{z}_1 +
(i/2)\partial\overline{\partial}\varphi, \\
$$
$$
\Omega_2 = Re(e^{i\theta}dz_1\wedge d z_2),\hspace{.1in} \Omega_3 =
Im(e^{i\theta}dz_1\wedge d z_2),
$$
where $\theta$ is a real constant and $\varphi$ is a smooth function
on $M$ such that
\begin{equation}\label{kod-tr}
4i(Im(dz_1\wedge d\overline{z}_2) +iRe(z_1)(dz_1\wedge
d\overline{z}_1))\wedge\partial\overline{\partial}\varphi =
\partial\overline{\partial}\varphi\wedge\partial\overline{\partial}\varphi
\end{equation}
Conversely, any three forms $\Omega_1,\Omega_2,\Omega_3$ of the form
given above determine a para-hyperk\"ahler structure on $M$.
Moreover, its metric  is flat if and only if $\varphi$ is constant.

 \end{thm}

 Note that any primary Kodaira surface is a toric bundle over an elliptic curve and the
pull-back of any smooth function on the base curve gives a solution
to (\ref{kod-tr}). This shows that the moduli space of
para-hyperk\"ahler structures on a primary Kodaira surface is
infinite dimensional, which is in sharp contrast with the positive
definite case.

\smallskip

Non-compact examples of para-hyperk\"ahler mnifolds can be
constructed   by means of the so-called Walker manifolds.

Recall that a Walker manifold is a triple $(M,g,\mathcal{D})$, where
$M$ is a smooth manifold, $g$ an indefinite metric,  and
$\mathcal{D}$ a parallel null distribution. The local structure of
such manifolds has been described by A.Walker \cite{Walker1950a} and
we refer to \cite{Derdzinski} for a coordinate-free version of his
theorem. Of special interest are the Walker metrics on $4$ manifolds
for which $\mathcal{D}$ is of dimension $2$ since they appear in
several specific pseudo-Riemannian structures. For example, the
metric of every para-hyperk\"ahler structure is Walker,
$\mathcal{D}$ being the $(+1)$-eigenbundle of either of its product
structures.

 According to \cite{Walker1950a}, for every Walker $4$-manifold $(M,g,\mathcal{D})$ with $dim\,\mathcal{D}=2$, there
exist local coordinates $(x,y,z,t)$ around any point of $M$ such
that the matrix of $g$ has the form
\begin{equation}\label{Walker metric}
g_{(x,y,z,t)} =\left(\begin{array}{cccc}
    0&0&1&0\\
    0&0&0&1\\
    1&0&a&c\\
    0&1&c&b\\
\end{array}\right)
\end{equation}
for some smooth functions $a$, $b$ and $c$. Then a local orthonormal
frame of $TM$ can be defined by
\[
\begin{array}{lll}
{{\bf e}_1=\displaystyle{\frac{1-a}{2}}
\partial_x + \partial_z}, \quad\quad\,\,\, {\bf
e}_2=\displaystyle{\frac{1-b}{2}} \partial_y + \partial_t
-c \partial_x, \\
                                                 \\
{{\bf e}_3=-\displaystyle{\frac{1+a}{2}} \partial_x  + \partial_z, \quad\quad {\bf
e}_4=-\displaystyle{\frac{1+b}{2}}
\partial_y
          +\partial_t - c \partial_x}.
\end{array}
\]
Let $\{J_1,J_2,J_3\}$ be the (local) almost para-hypercomplex
structure for which $J_1{\bf e_1}={\bf e_2}, J_1{\bf e_3}={\bf
e_4}$, $J_2{\bf e_1}={\bf e_3}, J_2{\bf e_2}=-{\bf e_4}$, $J_3{\bf
e_1}={\bf e_4}, J_3{\bf e_2}={\bf e_3}$. This structure is
compatible with the Walker metric $g$, thus we have an almost
para-hyperhermitian structure, called proper in \cite{Mat}.

\smallskip

The next two results have been proved in \cite{DDGMMV2}.

\begin{thm}\label {PC} The structure
$(g,J_1,J_2,J_3)$ is para-hyperhermitian if and only if the functions $a$, $b$
and $c$ have the form
%\begin{equation}\label{eq:hyper-parahermitian}
\[
\begin{array}{l}
a=x^2K+x P+\xi, \\[0.1in]
b=y^2K+y T+ \eta, \\[0.1in]
c=xyK+\frac{1}{2}x T+\frac{1}{2}y P+ \gamma,
\end{array}
\]
%\end{equation}
where the capital and Greek letters stand for arbitrary smooth functions of
$(z,t)$.
\end{thm}

\begin{thm}\label{HK}
The structure $(g,J_1,J_2,J_3)$ is para-hyperk\"ahler if and only if the
functions $a$, $b$ and $c$ do not depend on $x$ and $y$.
\end{thm}

 In particular, the above theorem shows that the neutral K\"ahler metrics considered by Petean \cite{P} are all
para-hyperk\"ahler and hence self-dual and Ricci-flat.

 By the Kamada results mentioned above, any compact para-hyperk\"ahler surface admits
two parallel, null and orthogonal vector fields. Conversely.
Theorem~\ref{HK} together with the Petean's classification of
neutral Ricci-flat K\"ahler surfaces (\cite{P}) leads to the
following
\begin{thm}
Let $(M,g)$ be a compact oriented neutral 4-manifold with two
parallel , null and orthogonal vector fields. Then $M$ admits a
para-hyperk\"ahler structure $\{g,J_1,J_2,J_3\}$, so $(M,J_1)$ is
biholomorphic to a complex torus or a primary Kodaira surface.
\end{thm}
A detailed proof of this theorem will appear in \cite{inprogress}.

\section{Existence of para-hyperhermitian metrics}

In this section we discuss the problem of existence of a metric
compatible with a given (almost) para-hypercomplex structure.

A (linear) para-hypercomplex structure on a vector space $V$ is a
triple $\{J_1,J_2,J_3\}$ of endomorphisms of $V$ satisfying the
relations (\ref{3.2}). Note that such structures exist on any
even-dimensional vector space. A metric $g$ on $V$ is called
compatible with the structure $\{J_1,J_2,J_3\}$ if the identities
(\ref{3.3}) are satisfied. In this case we say that
$\{g,J_1,J_2,J_3\}$ is a (linear) para- hyperhermitian structure on
$V$.

  If we are given a hypercomplex structure $\{J_1,J_2,J_3\}$  and $g$ is any
positive definite metric on $V$, then
$h(X,Y)=g(X,Y)+g(J_1X,J_1Y)+g(J_2X,J_2Y)+g(J_3X,J_3Y)$ is a positive
definite metric compatible with  $\{J_1,J_2,J_3\}$. In the case of a
para-hypercomplex structure, some authors suggest, by an analogy, to
consider the bilinear form
$h(X,Y)=g(X,Y)+g(J_1X,J_1Y)-g(J_2X,J_2Y)-g(J_3X,J_3Y)$ where $g$ is
a metric. This symmetric form is compatible with the given
para-hypercomplex structure but it may be degenerate.

\smallskip

\noindent {\bf Example}.  Let $e_1,e_2,e_3,e_4$ be the standard
bases of ${\Bbb R}^4$ and let $\{J_1,J_2,J_3\}$ be the
para-hypercomplex structure on  ${\Bbb R}^4$ for which $J_1e_1=e_2,
J_1e_3=e_4$, $J_2e_1=e_3, J_2e_2=-e_4$, $J_3e_1=e_4, J_3e_2=e_3$. If
$g$ is the standard metric on ${\Bbb R}^4$, then the endomorphisms
$J_1,J_2,J_3$ are isometries of $g$, hence the form
$h(X,Y)=g(X,Y)+g(J_1X,J_1Y)-g(J_2X,J_2Y)-g(J_3X,J_3Y)$ is
identically zero. Similarly, if $g$ is the metric for which
$e_1,...,e_4$ is an orthogonal basis with $g(e_1,e_1)=g(e_3,e_3)=1$,
$g(e_2,e_2)=g(e_4,e_4)=-1$ (in this case $J_1$ and $J_3$ are
anti-isometrices, while $J_2$ is an isometry).

\smallskip

  The next observation is implicitly contained in \cite{Kamada}.

\begin{lm}
Let $\{J_1,J_2,J_3\}$ be a para-hypercomplex structure on a vector
space $V$. Let $V^{\pm}$ be the $\pm 1$-eigenspace of the
endomorphism $J_2$. Then there is a bijective correspondence between
the set of non-degenerate skew-symmetric $2$-forms on the space
$V^{\pm}$ and the set of metrics on $V$ compatible with the given
para-hypercomplex structure.
\end{lm}

{\it Proof}. Let $h$ be a non-degenerate skew-symmetric $2$-form on
$V^+$. Extend this form to a form on the whole space $V$ setting
$h(V,V^-)=h(V^-,V)=0$. Now set $g(X,Y)=h(X,J_1Y)+h(Y,J_1X)$ for
$X,Y\in V$. Then $g$ is a symmetric bilinear form on $V$ and
$g(J_1X,J_1Y)=g(X,Y)$. Note that the spaces $V^{\pm}$ are
$g$-isotropic since $J_1$ interchanges $V^+$ and $V^-$. Let $X=X^+ +
X^-$, $Y=Y^+ +Y^-$ be the $V^{\pm}$-decomposition of arbitrary
vectors $X,Y\in V$. Then
$g(J_2X,J_2Y)=h(J_2X,J_3Y)+h(J_2Y,J_3X)=h(X^+,J_3Y^-)+h(Y^+,J_3X^-)$
since $J_3$ interchanges $V^+$ and $V^-$. On the other hand,
$g(X,Y)=h(X^+,J_1Y^-)+h(Y^+,J_1X^-)=-h(X^+,J_1J_2Y^-)-h(Y^+,J_1J_2X^-)
=-h(X^+,J_3Y^-)-h(Y^+,J_3X^-)$. Thus $g(J_2X,J_2Y)=-g(X,Y)$. It
follows that $g(J_3X,J_3Y)=-g(X,Y)$ since $J_3=J_1J_2$. Finally, the
identity $g(X,Y)=h(X^+,J_1Y^-)+h(Y^+,J_1X^-)$ and the fact that $h$
is non-degenerate on $V^+$ imply that $g$ is non-degenerate.

  Conversely, let $g$ be a metric on $V$ compatible with the para-hypercomplex structure
$\{J_1,J_2,J_3\}$. Then the spaces $V^{\pm}$ are $g$-isotropic. It
follows that $h(A,B)=\frac{1}{2}g(J_1A,B)$, $A,B\in V^+$, is a
non-degenerate skew-symmetric $2$-form. It is easy to check that $h$
yields the metric $g$.

\smallskip

  The proof above gives also the following

\begin{prop}\label{com}
Let $\{J_1,J_2,J_3\}$ be an almost para-hypercomplex structure on a
four-manifold $M$. Let $V^{\pm}$ be the subbundle of $TM$
corresponding to the eigenvalue $\pm 1$ of $J_2$. The manifold $M$
admits a metric $g$ compatible with the given para-hypercomplex
structure if and only if the bundle $V^{\pm}$ is orientable.
\end{prop}

The bundle $V^{\pm}$ is orientable iff the linear bundle $\Lambda^2
V^{\pm}$ is trivial. It follows that if $H^1(M,{\cal C}^{\ast})=0$
where ${\cal C}^{\ast}$ is the sheaf of non-vanishing smooth
real-valued function on $M$, then $V^{\pm}$ is orientable, hence
$\{J_1,J_2,J_3\}$ admits a compatible metric.

\smallskip

It is well-known that for every vector bundle there is a double
cover of its base such that the pull-back bundle is orientable.
Therefore we have the following

\begin{cor}\label{cover}
For any almost para-hypercomplex structure on a four-manifold $M$
there is a double cover of $M$ such that the pull-back
para-hypercomplex structure on it admits a compatible metric.
\end{cor}

Proposition~\ref{com} and the fact that any bundle on a simply
connected manifold is orientable imply

\begin{cor}
Any almost para-hypercomplex structure on a simply connected
four-manifold $M$ admits a compatible metric.
\end{cor}

We should emphasize that, in contrast to the definite case, not
every para-hypercomplex structure admits a compatible metric.
Examples of such structures on Inoue surfaces of type $S^{-}$ will
be provided in the last section. Other examples can be constructed
on hyperelliptic surfaces \cite{inprogress}.

\smallskip

  The next fact is well-known \cite{Bl-V} and easy to prove.

\begin{lm}
Let $\{g, J_1,J_2,J_3\}$ be a para-hyperhermitian structure on a
vector space $V$. A vector $w\in V$ is $g$-non-isotropic if and only
if $w,J_1w,J_2w,J_3w$ is a basis of $V$.
\end{lm}

{\it Proof}. The vectors $w,J_1w,J_2w,J_3w$ are $g$-orthogonal.
Thus, if $w$ is non-isotropic, they form a basis. Conversely,
suppose that $w,J_1w,J_2w,J_3w$ is a bases. Take a vector $e_1\in V$
with $||e_1||_g=1$. Then $e_1,e_2=J_1e_1, e_3=J_2e_1, e_4=J_3e_1$ is
a $g$-orthonormal basis of $V$. Let $(w_1,w_2,w_3,w_4)$ be the
coordinates of $w$ with respect to this basis. Then the coordinates
of $J_1w,J_2w,J_3w$ are $J_1w=(-w_2,w_1,-w_4,w_3)$,
$J_2w=(w_3,-w_4,w_1,-w_3)$, $J_3w=(w_4,w_3,w_2,w_1)$. It follows
that the transition matrix from the bases $(w,J_1w,J_2w,J_3w)$ to
the bases $(e_1,e_2,e_3,e_4)$ has determinant equal to
$(w_1^2+w_2^2-w_3^2-w_4^2)^2=||w||_g^4$. Hence $w$ is non-isotropic.

\smallskip

\smallskip

  This observation implies that any metric compatible with a
para-hyper-complex structure is of split signature.

\begin{lm}\label{confm}
Let $\{J_1,J_2,J_3\}$ be a para-hypercomplex structure on a vector
space $V$. Let $g$ and $h$ be two compatible metrics. If $w$ is an
$h$-non-isotropic vector, then it is also $g$-non-isotropic and
$g=\lambda\,h$, where $\lambda=g(w,w)/h(w,w)$.
\end{lm}

{\it Proof}. It is clear that the identity
$g(X,Y)=g(w,w)/h(w,w)\,h(X,Y)$ holds when $X=Y=w,J_1w,J_2w,J_3w$.
Hence it holds for every $X,Y\in V$ since $w,J_1w,J_2w,J_3w$ is a
basis which is $g$- and $h$-orthogonal. In particular, $g(w,w)\neq
0$.

\begin{prop}\label{confclass}
If $\{J_1,J_2,J_3\}$ is an almost para-hypercomplex structure on a
four-manifold $M$ and $g$,$h$ are two compatible metrics, then there
exists a unique non-vanishing smooth function $f$ on $M$ such that
$g=f\,h$.
\end{prop}

{\it Proof}. Every point of $M$ has a neighbourhood  with an
$h$-non-isotropic vector field $W$ on it and the proposition follows
form Lemma~\ref{confm}

\smallskip

   Corrolary~\ref{cover} and Proposition~\ref{confclass} imply the following

\begin{prop}
Every para-hypercomplex structure on a $4$-manifold $M$ determines a
conformal class up to a double cover of $M$.
\end{prop}

\section{Inoue surfaces of type $S^{\pm}$}

   It is well-known \cite{Bl-V,kuche} that the Inoue surfaces of
type $S^+$ admit para-hyperhermitian structures. In this section, we
show that, in contrast, any Inoue surface of type $S^-$ has a
para-hypercomplex structure which does not admit a compatible
metric. Before that we recall the definition of the Inoue surfaces
of type $S^{\pm}$.

 Let $p,q,r $ be integers, $t$ a complex number, and $N\in SL(2,\mathbb{Z})$  a matrix with
eigenvalue $\alpha>1$ and $1/\alpha$. Denote by $\mathbb{H}$ the
upper half-plane of the complex plane $\mathbb{C}$,
The Inoue surface  $S^+_{p,q,r,t, N}$ is obtained as a quotient of $\mathbb{H}\times\mathbb{C}$ by the action of the group generated by following transformations:\\
$\phi_0(z,w)\rightarrow(\alpha z,w+t)$\\
$\phi_i(z,w)\rightarrow (z+a_i,w+b_iz+c_i)\,, \,\,\,\, i = 1, 2$\\
$\phi_3(z,w)\rightarrow(z,w+A)$,\\
where $(a_1,a_2)$ and $(b_1,b_2)$ are real eigenvectors of $N$
corresponding to $\alpha$ and $1/\alpha$, and $A=(b_1a_2-b_2a_1)/r$.
Here, the constants $c_i$ are real numbers determined by $a_i, b_i,
 p,q,r$,\,\,and  the eigenvalues of $N$.\\
The (1,0)-forms
$$\theta_1=\frac{dz}{Im \,z}\,\,\,\,\, \text{and} \,\,\,\,\, \theta_2 = dw-\frac{Im\, w-s\ln(Im \,z)}{Im \,z} dz$$
where $s= Im\, t/\ln \alpha$ are invariant under this action and the
corresponding
dual (1,0)-vector fields are:\\
$$E_1=(Im\,z)\frac{\partial}{\partial z}+(Im\, w - s \ln(Im\,z))
\frac{\partial}{\partial w} \,\,\,\,\,\, \text{ and} \,\,\,\,\,\,
E_2=\frac{\partial}{\partial w}.$$ It is easy to see that
$$d\theta_1 = (-1/2i)\theta_1\wedge \overline{\theta_1},\,\,\,\,\,\,\,\, d\theta_2= (1/2i)(\theta_1\wedge \theta_2-\theta_1\wedge \overline{\theta_2} + s\theta_1 \wedge \overline{\theta_1}).$$
From here one gets
$$d(\theta_1\wedge\theta_2)=
-(Im\,\theta_1)\wedge\theta_1\wedge\theta_2,$$ thus the $(2,0)$-form
$\Omega = \theta_1\wedge\theta_2$ satisfies the relation
$d\Omega=-Im\,\theta_1\wedge\Omega$.\\
%This latter can also be seen by using that
%$$\Omega=\frac{dz\wedge dw}{Im z}\,\,\,\,\, \text{and} \,\,\,\,\,Im\,\theta_1=d(\ln Im z).$$

 Set $\Omega_1=Re(\theta_1\wedge\overline\theta_2)$. Then one can
 check that $$d\Omega_1=-Im\,\theta_1\wedge\Omega_1,\quad
 \Omega_1^2=-(Re\,\Omega)^2=-(Im\,\Omega)^2=\frac{1}{2}\theta_1\wedge\overline\theta_1
 \wedge\theta_2\wedge\overline\theta_2.$$
Therefore the triple $(\Omega_1,Re\,\Omega,Im\,\Omega)$ defines a
para-hyperhermitian structure on $S^+_{p,q,r,t, N}$.

 The definition of Inoue surfaces of type $S^-$ is the same as those of type $S^+$, but in this case $\phi_0$ is
defined as $\phi_0(z,w)\rightarrow(\alpha z,-w)$. It is clear that
any surface $S^-$ is a quotient of a certain surface $S^+$ with
$t=0$ by the action of the involution $\sigma: S^+\rightarrow S^+
\,\,\,\,\text{given by} \,\,\,\,\sigma(z,w) = (z, -w)$. Then for the
$(1,0)$-forms $\theta_1$ and $\theta_2$ defined above, we have
$\sigma^{\ast}\theta_1=\theta_1,
\,\,\,\,\,\sigma^{\ast}\theta_2=-\theta_2$. Therefore
$\sigma^*\Omega_1=-\Omega_1$ and $\sigma^*\Omega=-\Omega$, hence the
para-hyperhermitian structure on $S^+$ defined above does not
descend to  $S^-$. Nevertheless, we show below that the surface
$S^-$ admits a para- hypercomplex structure with no compatible
metric. Notice first that the para-hypercomplex structure on $S^+$
defined by the para-hyperhermitian structure
$(\Omega_1,Re\,\Omega,Im\,\Omega)$ does descend to $S^-$. Indeed,
the map $\sigma$ is an anti-isometry with respect to the metric of
this structure since it preserves the complex structure and
$\sigma^*(\Omega_1)=-\Omega_1$. This, the identity
$\sigma^*\Omega=-\Omega$ and the fact that $\sigma$ is an involution
imply that $\sigma$ preserves also the two product structures. Hence
the para-hypercomplex structure descends to $S^-$.  On the other
hand, if we suppose that there is a  metric on $S^-$ compatible with
the induced para- hypercomplex structure, its pull-back would be a
metric compatible with the para- hypercomplex structure  on $S^+$.
But, according to Proposition~\ref{confclass}, in real dimension 4,
any two metrics compatible with the same para-hypercomplex structure
are conformally equivalent. So there would be a nowhere vanishing
real-valued function $f$ on $S^+$ for which $f\Omega$ is the
pull-back of the $(2,0)$ - form on $S^-$ associated with the
para-hyperhermitian structure there. Since $f\Omega$ is
$\sigma$-invariant and $\sigma^*(\Omega)=-\Omega$, we have $f(\sigma
(x))=-f(x)$ for every $x\in S^+$. But this contradicts to the fact
that $f$ has a fixed sign. So the para-hypercomplex structure on
$S^-$ defined above does not admit a compatible metric.

 \vskip 20pt

\noindent Johann Davidov and Oleg Mushkarov

\noindent Institute of Mathematics and Informatics

\noindent Bulgarian Academy of Sciences

\noindent 1113 Sofia, Bulgaria

\noindent jtd@math.bas.bg,  muskarov@math.bas.bg

\medskip

\noindent Gueo Grantcharov and Miroslav Yotov

\noindent Department of Mathematics

\noindent Florida International University

\noindent Miami, FL 33199

\noindent grantchg@fiu.edu, yotovm@fiu.edu

\end{document}